\documentclass[11pt, reqno]{amsart}

\usepackage{amsmath}
\usepackage{fullpage}
\usepackage{amssymb}
\usepackage{amsthm}
\usepackage{bbm}
\usepackage{subfigure}
\usepackage{float}
\usepackage{graphicx}
\usepackage{color}
\usepackage{epstopdf}

\newtheorem{lemma}{Lemma}[section]
\newtheorem{theorem}[lemma]{Theorem}
\newtheorem{proposition}[lemma]{Proposition}

\newtheorem{definition}[lemma]{Definition}
\graphicspath{{../figs/}}

\DeclareMathOperator{\tr}{Tr}

\title{Low regularity full error estimates for the cubic nonlinear Schr\"odinger equation}

\author{Lun Ji}
\address{Department of Mathematics, Universit\"{a}t Innsbruck, Technikerstr.~13, 6020 Innsbruck, Austria (L.~Ji)}
\email{lun.ji@uibk.ac.at}

\author{Alexander Ostermann}
\address{Department of Mathematics, Universit\"{a}t Innsbruck, Technikerstr.~13, 6020 Innsbruck, Austria (A.~Ostermann)}
\email{alexander.ostermann@uibk.ac.at}

\author{Fr\'ed\'eric Rousset}
\address{Universit\'e Paris-Saclay, CNRS, Laboratoire de Math\'ematiques d'Orsay (UMR 8628), 91405 Orsay Cedex, France (F.~Rousset)}
\email{frederic.rousset@universite-paris-saclay.fr}

\author{Katharina Schratz}
\address{LJLL (UMR 7598), Sorbonne Universit\'e, UPMC, 4 place Jussieu, 75005, Paris, France (K.~Schratz)}
\email{katharina.schratz@ljll.math.upmc.fr}

\begin{document}

\begin{abstract}
For the numerical solution of the cubic nonlinear Schr\"{o}dinger equation with periodic boundary conditions, a pseudospectral method in space combined with a filtered Lie splitting scheme in time is considered. This scheme is shown to converge even for initial data with very low regularity. In particular, for data in $H^s(\mathbb T^2)$, where $s>0$, convergence of order $\mathcal O(\tau^{s/2}+N^{-s})$ is proved in $L^2$. Here $\tau$ denotes the time step size and $N$ the number of Fourier modes considered. The proof of this result is carried out in an abstract framework of discrete Bourgain spaces, the final convergence result, however, is given in $L^2$. The stated convergence behavior is illustrated by several numerical examples.
\end{abstract}

\maketitle

\section{Introduction}

The  periodic cubic nonlinear Schr\"{o}dinger equation (NLS)
\begin{equation}\label{nls-dd}
iu_t=-\Delta u-\mu|u|^2u,\quad (t,x)\in\mathbb{R}\times\mathbb{T}^d,\qquad \mu=\pm1
\end{equation}
on the $d$-dimensional torus $\mathbb T^d$ is typically integrated numerically using Strang splitting in time and some pseudospectral method (usually FFT) in space~\cite{Faou}. This approach, however, requires a certain amount of smoothness of the solution \cite{Lubich} and becomes ineffective for non-smooth initial data. In the last decade, there has been a significant effort to develop new numerical integration schemes that demand less regularity. In particular, various filtered exponential integrators and splitting methods were constructed that treat different frequencies in the solution differently \cite{Bao, LiWu, OstS, OstY}. Most of these methods have been considered in a standard Sobolev space framework, where the Sobolev index $s$ has to satisfy the restriction $s>d/2$ due to the bilinear estimate
$$
\|fg\|_{H^s} \le c_{s,d} \|f\|_{H^s}\|g\|_{H^s}, \quad f,g\in H^s(\mathbb T^d)
$$
required in the analysis.

More recently, with the help of discrete Bourgain techniques, much less regular initial data can be considered~\cite{Ji, JiO, Ost, Ost1}. For instance, in dimensions $d=1$ and $d=2$, initial data in $H^s(\mathbb T^d)$ are admissible for $s>0$ and lead to the temporal order $s/2$ in $L^2$. The author of \cite{Wu} improved this convergence rate by using more sophisticated filtering techniques. The resulting algorithm treats the low and high frequencies separately and also preserves the mass. However, in the papers above, only the time discretization was considered. Note that in the Sobolev setting, some convergence results for full discretizations are available in the literature~\cite{Bao, Ignat, LiWu, OstY}. However, nothing was known for very low regularity data in the Bourgain space setting.

The aim of this paper is to fill this gap. For notational simplicity, we will fix the dimension to $d=2$. Note, however, that our results can easily be extended to other dimensions. For the time discretization, this was already done in \cite{JiO}. Henceforth, we consider the periodic cubic nonlinear Schr\"{o}dinger equation
\begin{equation}\label{nls}
iu_t=-\Delta u-\mu|u|^2u,\quad (t,x)\in\mathbb{R}\times\mathbb{T}^2,\qquad \mu=\pm1,
\end{equation}
defined on the 2-dimensional torus $\mathbb T^2$. Both, the focusing case ($\mu=1$) and the defocusing case ($\mu=-1$) can be studied.

Together with this problem, we also consider the projected NLS
\begin{equation}\label{projeq}
i\partial_t u^\theta =-\Delta u^\theta -\mu\Pi_\theta(|\Pi_\theta  u^\theta |^2\Pi_\theta  u^\theta ), \quad u^\theta(0)=\Pi_\theta  u(0)
\end{equation}
for some appropriate parameter $\theta>0$. Here, the projection operator $\Pi_\theta $ for $\theta>0$ is defined by the Fourier multiplier:
\begin{equation}\label{proj}
\Pi_\theta =\overline{\Pi}_\theta=\chi\left(\dfrac{-i\nabla}{\theta^{-\frac12}}\right),
\end{equation}
where $\chi$ is the characteristic function of the square $[-1,1)^2$.

Moreover, we shall need the following projection operator which is related to the discrete Fourier transform (DFT).

\begin{definition}\label{tndefine}
For every even number $N>0$, $T_N$ denotes the following operator, acting on continuous functions $u: \mathbb{T}^2 \to \mathbb C$:
\begin{equation}\label{tndef}
T_N(u)(x_1,x_2)=\frac1{N^2}\sum\limits_{k_1=-\frac N2}^{\frac N2-1}\sum\limits_{k_2=-\frac N2}^{\frac N2-1}\mathcal{F}_N(u)(k_1,k_2)e^{ik_1x_1+ik_2x_2},
\end{equation}
where $\mathcal{F}_N(u)$ is the two-dimensional discrete Fourier transform, i.e.
\begin{equation}\label{dft}
\mathcal{F}_N(u)(k_1,k_2)=\sum\limits_{j_1=-\frac N2}^{\frac N2-1}\sum\limits_{j_2=-\frac N2}^{\frac N2-1}u(\tfrac{2\pi j_1}{N},\tfrac{2\pi j_2}{N})e^{-\tfrac{2i\pi k_1j_1}{N}-\tfrac{2i\pi k_2j_2}{N}}, \quad -\tfrac N2\leq k_1,k_2\leq\tfrac N2-1.
\end{equation}
\end{definition}
Note that the result, $T_N(u)$, is again a continuous function on the torus $\mathbb{T}^2$. 

For the space discretization of \eqref{projeq} a pseudospectral method is considered. Specifically, we discretize the torus $\mathbb T^2$ in both directions by an equidistant grid and use the DFT to map between physical and frequency space. This means in particular that the operator $T_N$ has to be applied to the nonlinearity. With Definition~\ref{tndefine} at hand, we can now define a fully discrete filtered Lie splitting scheme
\begin{equation}\label{liespace}
u_{n+1}=\Psi^\tau_N(u_n)=e^{i\tau\Delta}\Pi_\theta T_N(e^{\mu i\tau|\Pi_\theta  u_n|^2}\Pi_\theta  u_n),\quad u_0=\Pi_\theta  u(0),\qquad \theta=\max(\tau,4N^{-2})
\end{equation}
that will be analyzed in this paper. Note that, by construction, $u_n=\Pi_\theta u_n$. By the use of $T_N$, the error caused by aliasing is also included and will be analyzed. This goes beyond our previous work in \cite{Ji}. 

In this paper, we will prove the following result:

\begin{theorem}\label{mainthm}
For $s_0\in(0,2]$ and initial data $u_0\in H^{s_0}(\mathbb{T}^2)$, let $u$ be the exact solution of \eqref{nls} in the Bourgain space $X^{s_0,b_0}(T)$ with initial data $u_0$ for appropriate $T>0$ and $b_0>\tfrac12$ (see also \cite{Ji}). Further, let $u_n$ denote the numerical solution defined by the scheme \eqref{liespace}. Then, we have the following error estimate: there exist $\tau_0>0$, $N_0\in\mathbb{N}$ and $C_T>0$ such that for every time step size $\tau\in(0,\tau_0]$ and (even) number of Fourier modes $N>N_0$ per dimension, the global error satisfies
$$
\Vert u_n-u(t_n)\Vert_{L^2(\mathbb{T}^2)}\leq C_T\,\theta^{\frac{s_0}2},\quad 0\leq n\tau\leq T,
$$
where the constant $C_T$ depends on $\tau_0$, $N_0$ and $T$, but is independent of $n$, $\tau$ and $N$.
\end{theorem}

For notational simplicity, we have chosen the same number of Fourier modes in each direction. This number equals to the number of spatial grid points in each dimension because of the use of DFT. To further simplify our notations, we always take the defocusing case ($\mu=-1$) in the following. We stress, however, that our analysis given below remains true for $\mu=1$ as long as the solution exists on $[0, T]$.

\subsection*{Outline of the paper.} The paper is organized as follows. In section \ref{sectionanalysis}, we briefly recall the main steps of the analysis of the Cauchy problem for (\ref{nls}), and we prove an estimate on the difference between the exact solutions of (\ref{nls}) and (\ref{projeq}). In section \ref{sectiondiscbourg}, we recall the main properties of the discrete Bourgain spaces $X^{s,b}_{\tau}$. In section \ref{sectionlocal}, we analyze the local spatial error of the scheme (\ref{liespace}); in section~\ref{sectionglobal}, we give global error estimates and prove our main result, Theorem \ref{mainthm}. Numerical examples, which are given at the end of the paper, illustrate our convergence result.

\subsection*{Notations.}
We denote the time step size by $\tau$, the number of the spatial grid points by $N$ (in any of the two spatial directions), and we set $\theta=\max(\tau,4N^{-2})$. We also tacitly assume that the parameters $s_0,~b_0,~b_1$ satisfy
\begin{equation}\label{s0b0}
s_0>0,\quad b_0\in\big(\tfrac12,\max(\tfrac34,\tfrac12+\tfrac14s_0)\big),\quad b_1=1-b_0,
\end{equation}
and we use the Japanese bracket $\langle\,\cdot\,\rangle=(1+|\cdot|^2)^\frac12$ notation. 

The estimate $A\lesssim B$ means $A\leq CB$, where $C$ is a generic constant; in particular, $C$ is independent of $\tau \in (0, 1]$ and $N$. The symbol $\lesssim_\gamma$ emphasizes that the constant $C$ depends in particular on $\gamma$. Moreover, $A\sim B$ means that $A\lesssim B\lesssim A$.

Finally, for sequences $(u_n)_{n\in\mathbb{Z}^m} \in X^{\mathbb{Z}^m}$, where $X$ is a Banach space with norm $\|\cdot \|_X$, we employ the usual norms
\begin{equation}\label{lpnorm}
\|u_n\|_{l^p_\zeta X}=\Big(\zeta^m\sum_{n\in\mathbb{Z}^m}\|u_n\|_X^p\Big)^\frac1p,\quad \|u_n\|_{l^\infty_\zeta X}=\sup_{n\in\mathbb{Z}^m}\|u_n\|_X.
\end{equation}

\section{Analysis of the exact solution}\label{sectionanalysis}

In this section, we will discuss the Cauchy problem for (\ref{nls}) at low regularity. The use of Bourgain spaces and some subtle multilinear estimates is crucial for that. We shall then also discuss the Cauchy problem for the projected equation \eqref{projeq} and an estimate between the solution $u$ of \eqref{nls} and the solution $u^\theta $ of \eqref{projeq}.

Let us recall the definition of Bourgain spaces first. For a function $u(t,x)$ on $\mathbb{R}\times\mathbb{T}^2$, $\tilde{u}(\sigma,k)$ stands for its time-space Fourier transform, i.e.
$$
\tilde{u}(\sigma,k)=\int_{\mathbb{R}\times\mathbb{T}^2}u(t,x)e^{-i\sigma t-i\langle k,x\rangle}dx dt,
$$
where $\langle\cdot,\cdot\rangle$ denotes the inner product in $\mathbb R^2$. The inverse transform is given by
$$
u(t,x)=\frac1{4\pi^2}\sum\limits_{k\in\mathbb{Z}^2}\hat{u}_k(t)e^{i\langle k,x\rangle}
$$
with the Fourier coefficients $\hat{u}_k(t)=\frac1{2\pi}\int_{\mathbb{R}}\tilde{u}(\sigma,k)e^{i\sigma t}d\sigma$.

This way, we can define the Bourgain space $X^{s,b}= X^{s,b}(\mathbb{R}\times\mathbb{T}^2)$, which is a Banach space equipped with norm
$$
\Vert u\Vert_{X^{s,b}}=\Vert\langle k\rangle^s\langle\sigma + |k|^2\rangle^b\tilde{u}(\sigma,k)\Vert_{L^2 l^2}.
$$

We also recall the definition of a localized version of this space: $X^{s,b}(I)$, where $I\subset\mathbb{R}$ by
\begin{equation}\label{bourgainloc}
\Vert u\Vert_{X^{s,b}(I)}=\inf\{\Vert v\Vert_{X^{s,b}}~\big|~v|_I=u\},
\end{equation}
and we write $X^{s,b}(T)=X^{s,b}([0, T])$ for short.

Then, let us recall some well-known properties of Bourgain spaces.

\begin{lemma}\label{contiprop}
For any $s\in\mathbb{R}$, $b>\tfrac12$, $s_0,~b_1$ fulfilling the conditions in \eqref{s0b0} and any permutation $\gamma$ of $\{1, 2, 3\}$, we have that
\begin{align}
\label{d}\Vert u\Vert_{L^\infty H^s}&\lesssim_b\Vert u\Vert_{X^{s,b}},\\
\label{eqmult}\Vert u\overline{v}w\Vert_{X^{s_0,-b_1}}&\lesssim_{s_0,b_1}\Vert u\Vert_{X^{s_0,b_1}}\Vert v\Vert_{X^{s_0,b_1}}\Vert w\Vert_{X^{s_0,b_1}},\\
\label{eqmult2}\Vert u_1\overline{u}_2u_3\Vert_{X^{0,-b_1}}&\lesssim_{s_0,b_1} \Vert u_{\gamma(1)}\Vert_{X^{s_0,b_1}}\Vert u_{\gamma(2)}\Vert_{X^{s_0,b_1}}\Vert u_{\gamma(3)}\Vert_{X^{0,b_1}}.
\end{align}
\end{lemma}

These estimates are proven, for example, in \cite[sect.~2.6]{Tao} and in \cite{Burq,Ji}.

We shall now establish the existence and uniqueness of the exact solution of \eqref{nls}.

\begin{theorem}\label{theoexist}
For $s_0,~b_0$ fulfilling the conditions in \eqref{s0b0}, let $u_0\in H^{s_0}(\mathbb{T}^2)$. Then, there exists $T^*\in (0, +\infty]$ and a unique maximal solution $u\in  X^{s_0,b_0}(T) \subset \mathcal{C}([0,T],H^{s_0}(\mathbb{T}^2))$ for every $T \in (0, T^*)$ satisfying \eqref{nls}.
\end{theorem}

We refer to \cite[Theorem 2.3]{Ji} for the proof of this theorem.

To close this section, we will give another useful theorem. For its proof, we refer to \cite[sect.~2]{Ji}.

\begin{theorem}\label{theou-utau}
Let $u_0\in H^{s_0}$ and $u\in X^{s_0,b_0}(T)$ for any $T\in(0,T^*)$ be the solution of \eqref{nls} given by Theorem \ref{theoexist}. Then there exists $\theta_0>0$ such that for every $\theta\in(0,\theta_0]$, there is a unique solution $u^\theta$ of \eqref{projeq} which is also defined on $[0,T]$ satisfying $u^\theta\in X^{s_0,b_0}(T)$. Moreover, we also have that
\begin{equation} \label{diffutau}
\Vert u-u^\theta\Vert_{X^{0,b_0}(T)}\lesssim C_T\theta^\frac{s_0}2, \qquad \Vert u^\theta\Vert_{X^{s_0, b_0}(T)}\leq C_T
\end{equation}
for some $C_T$ independent of $\theta$ for $\theta\in (0,\theta_0].$
\end{theorem}

In the following, thanks to the definition of the local Bourgain spaces (see \eqref{bourgainloc}), we shall still denote by $u^\theta $ an extension of the solution $u^\theta $ of \eqref{projeq} on $[0, T]$ such that
\begin{equation}\label{bmieux}
\|u^\theta \|_{X^{s_{0}, b_0}}\leq 2\|u^\theta\|_{X^{s_0, b_0}(T)} \lesssim C_T, \qquad
\|u^\theta\|_{X^{s_0,1}}\leq 2\|u^\theta\|_{X^{s_0,1}(T)}\lesssim C_{T,\sigma}\theta^{-\sigma}
\end{equation}
where $\sigma$ can be chosen arbitrarily close to $\frac12$, and the choice is independent of $s_0$. Note that $u^\theta $ is now defined globally, but it is a solution of \eqref{projeq} only on $[0, T]$.

\section{Discrete Bourgain spaces}\label{sectiondiscbourg}

We now recall the definition of the discrete Bourgain spaces, see also \cite{Ji,Ost}. We first take $(u_n(x))_n$ to be a sequence of functions on the torus $\mathbb{T}^2$, with its time-space Fourier transform
$$
\widetilde{u_n}(\sigma,k)=\tau\sum\limits_{m\in\mathbb{Z}}\widehat{u_m}(k)e^{im\tau\sigma},
$$
where
\begin{equation}\label{spaceft}
\widehat{u_m}(k)=\dfrac1{4\pi^2}\int_{\mathbb{T}^2}u_m(x)e^{-i\langle k,x\rangle}dx.
\end{equation}
In this framework, $\widetilde{u_{n}}$ is a $ 2 \pi/\tau$ periodic function in $\sigma$ and Parseval's identity reads
\begin{equation}\label{parseval}
\| \widetilde{u_{n}}\|_{L^2l^2}= \|u_{n}\|_{l^2_{\tau}L^2},
\end{equation}
where we use the shorthand notation
$$
\| \widetilde{u_{n}}\|_{L^2l^2}^2 = \int_{-\frac\pi\tau}^\frac\pi\tau \sum_{k \in \mathbb{Z}^2}
|\widetilde{u_{n}}(\sigma, k)|^2 \, d \sigma, \qquad
\|u_{n}\|_{l^2_{\tau}L^2}^2 = \tau \sum_{m \in \mathbb{Z}} \int_{\mathbb{T}^2}|u_{m}(x)|^2 \, dx.
$$

Then the discrete Bourgain space $X^{s,b}_\tau$ can be defined with two equivalent norms (see \cite{Ji} and \cite{Ost})
\begin{equation}\label{5}
\Vert u_n\Vert_{X^{s,b}_\tau}=\Vert\langle k\rangle^s\langle d_\tau(\sigma-|k|^2)\rangle^b\widetilde{u_n}(\sigma,k)\Vert_{L^2l^2}\sim\Vert\langle D_\tau\rangle^b\langle k\rangle^s (e^{-in\tau\Delta}u_n)\Vert_{l_\tau^2L^2},
\end{equation}
where $d_\tau(\sigma)=\frac{e^{i\tau\sigma}-1}{\tau}$, $(D_\tau(u_n))_n=\big(\frac{u_{n-1}-u_n}{\tau}\big)_n$. Note that for any fixed $(u_n)_n$, the norm is an increasing function for both $s$ and $b$.

We can also define the localized discrete Bourgain spaces, $X_\tau^{s,b}(I)$, with
$$
\Vert u_n\Vert_{X_\tau^{s,b}(I)}=\inf\{\Vert v_n\Vert_{X_\tau^{s,b}}~\big|~v_n =u_n, \, n \tau \in I\}.
$$

We directly get from the definition (\ref{5}) the elementary properties that  for any $s\geqslant s^\prime$ and $b\geqslant b^\prime$, we have
\begin{align}
\label{b}\Vert\Pi_\theta  u_n\Vert_{X^{s,b}_\tau}&\lesssim\tau^{b^\prime-b}\Vert\Pi_\theta  u_n\Vert_{X^{s,b^\prime}_\tau},\\
\label{t}\Vert\Pi_\theta  u_n\Vert_{X^{s,b}_\tau}&\lesssim \theta^{\frac{s^\prime-s}2}\Vert\Pi_\theta u_n\Vert_{X^{s^\prime,b}_\tau}.
\end{align}

Next, we will provide some further useful estimates.
\begin{lemma}\label{discprop}
For $s\in\mathbb{R},~\eta\in\mathcal{C}_c^\infty(\mathbb{R})$ and $\tau\in(0,1]$, we have that
\begin{align}
\label{9}\Vert\eta(\tfrac{n\tau}T)u_n\Vert_{X^{s,b^\prime}_\tau}&\lesssim_{\eta,b,b^\prime}T^{b-b^\prime}\Vert u_n\Vert_{X^{s,b}_\tau},\quad -\tfrac12<b\leq b^\prime<\tfrac12,~0<T=M\tau\leq1,~M\geq1,\\
\label{a}\sup\limits_{\delta\in[-4,4]}\Vert e^{i\tau\delta\Delta}u_n\Vert_{X^{s,b}_\tau}&\lesssim\Vert u_n\Vert_{X^{s,b}_\tau},\quad b\in\mathbb{R},\\
\label{y}\Vert u_n\Vert_{l_\tau^\infty H^s}&\lesssim_b\Vert u_n\Vert_{X^{s,b}_\tau},\quad b>\tfrac{1}{2},\\
\label{ydual}\Vert u_n\Vert_{X^{s,-b}_\tau}&\lesssim_{b}\Vert u_n\Vert_{l_\tau^1 H^s},\quad b>\tfrac{1}{2},\\
\label{q}\Vert U_n\Vert_{X^{s,b}_\tau}&\lesssim_{\eta,b}\Vert u_n\Vert_{X^{s,b-1}_\tau},\quad b>\tfrac{1}{2},
\end{align}
where
$$U_n(x)=\tau\eta(n\tau)\sum\limits_{m=0}^ne^{i(n-m)\tau\Delta}u_m(x).$$
\end{lemma}

We stress that all these estimates are uniform in $\tau$.
The proof of this lemma is also nearly the same as the one-dimensional case given in \cite[sect.~3]{Ost}. Therefore, we omit the details. Note that \eqref{ydual} is a dual version of \eqref{y}.

The key estimates for the analysis of the scheme are given in the following theorem.
\begin{theorem}\label{disckeythm}
For any $s_0,~b_1$ defined in \eqref{s0b0}, $\theta\geq\tau$, and $c>0$, we have
\begin{align}
\label{z}\Vert\Pi_{c\theta}u_n\Vert_{l^4_\tau L^4}&\lesssim_c\Vert u_n\Vert_{X_\tau^{\frac{s_0}2,b_1}},\\
\label{r}\Vert\Pi_{c\theta}(\Pi_{c\theta}u_{n,1}\Pi_{c\theta}\overline{u}_{n,2}\Pi_{c\theta}u_{n,3})\Vert_{X_\tau^{s_0,-b_1}}&\lesssim_c\Vert u_{n,1}\Vert_{X_\tau^{s_0,b_1}}\Vert u_{n,2}\Vert_{X_\tau^{s_0,b_1}}\Vert u_{n,3}\Vert_{X_\tau^{s_0,b_1}}, \\
\label{rbis}\Vert\Pi_{c\theta}(\Pi_{c\theta}u_{n,1}\Pi_{c\theta}\overline{u}_{n,2}\Pi_{c\theta}u_{n,3})\Vert_{X_\tau^{0,-b_1}}&\lesssim_c\Vert u_{n, \gamma(1)}\Vert_{X_\tau^{s_0,b_1}}\Vert u_{n,\gamma(2)}\Vert_{X_\tau^{s_0,b_1}}\Vert u_{n,\gamma(3)}\Vert_{X_\tau^{0,b_1}},
\end{align}
where $(u_n)_n$, $(u_{n,i})_n$, are functions in the corresponding spaces and $\gamma$ is any permutation of $\{1, 2, 3\}$.
\end{theorem}
The proof of \eqref{z} is given in \cite{JiO}. Moreover, we can see that \eqref{r} and \eqref{rbis} are the discrete counterparts of \eqref{eqmult} and \eqref{eqmult2}. The proof of these two estimates uses the same arguments as the corresponding proof for $c=1$, given in \cite{Ji}.

In particular, if we take $c=\frac19$, and substitute $\Pi_\theta u_{n,i}$, $i=1,2,3$, we have
\begin{align}
\label{r9}\Vert\Pi_\theta u_{n,1}\Pi_\theta\overline{u}_{n,2}\Pi_\theta u_{n,3}\Vert_{X_\tau^{s_0,-b_1}}&\lesssim\Vert u_{n,1}\Vert_{X_\tau^{s_0,b_1}}\Vert u_{n,2}\Vert_{X_\tau^{s_0,b_1}}\Vert u_{n,3}\Vert_{X_\tau^{s_0,b_1}},\\
\label{rbis9}\Vert\Pi_\theta u_{n,1}\Pi_\theta\overline{u}_{n,2}\Pi_\theta u_{n,3}\Vert_{X_\tau^{0,-b_1}}&\lesssim\Vert u_{n,\gamma(1)}\Vert_{X_\tau^{s_0,b_1}}\Vert u_{n,\gamma(2)}\Vert_{X_\tau^{s_0,b_1}}\Vert u_{n,\gamma(3)}\Vert_{X_\tau^{0,b_1}},
\end{align}
since $\Pi_{\frac\theta9}\Pi_\theta u_{n,i}=\Pi_\theta u_{n,i}$ and the frequency of $\Pi_\theta u_{n,1}\Pi_\theta\overline{u}_{n,2}\Pi_\theta u_{n,3}$ is at most $3\theta^{-\frac12}=(\frac\theta9)^{-\frac12}$.

\section{Local error estimates}\label{sectionlocal}

First of all, we can split the local error into temporal and spatial parts as
\begin{equation}\label{localee}
\begin{aligned}
\Psi^\tau_N(u^\theta(t_n))-u^\theta(t_{n+1})&=\big(\Psi^\tau(u^\theta(t_n))-u^\theta(t_{n+1})\big)+\big(\Psi^\tau_N(u^\theta(t_n))-\Psi^\tau(u^\theta(t_n))\big)\\
&=ie^{i\tau\Delta}\mathcal{E}_{tloc}(t_n,N,\tau,u^\theta )+ie^{i\tau\Delta}\mathcal{E}_{sloc}(t_n,N,\tau,u^\theta ),
\end{aligned}
\end{equation}
where $\Psi_N^\tau$ is defined in \eqref{liespace}, and $\Psi^\tau$ is defined as
\begin{equation}
\Psi^\tau(u_n)=e^{i\tau\Delta}\Pi_\theta(e^{-i\tau|\Pi_\theta  u_n|^2}\Pi_\theta  u_n),\quad u_0=\Pi_\theta  u(0),
\end{equation}
By \cite[Theorem~5.1]{Ji}, the temporal part has the estimate
\begin{equation}\label{timeloc}
\Vert\mathcal{E}_{tloc}(t_k,N,\tau,u^\theta)\Vert_{X^{0,-b_1}_\tau}\lesssim C_T\tau\theta^{\frac{s_0}2}.
\end{equation}

We recall that we denote in the following by $u^\theta$ the extension of the solution of \eqref{projeq} on $[0,T]$, see Theorem~\ref{theou-utau}. This will not make any difference on the error estimates since we only care about the error on $[0,T]$.

Before we estimate the local error, we shall first prove the boundedness of $u^\theta(t_n)$ in the norm of $X^{s,b}_\tau$ for suitable $s$ and $b$.

\begin{lemma}
Let $u^\theta  \in X^{s_0,b_0}$ be the extension of the solution of the projected NLS \eqref{projeq} given by Theorem \ref{theou-utau}, and define the sequence $u^\theta _n(x)=u^\theta(n\tau+t^\prime,x)$ for $t^\prime\in[0,4\tau]$. Then, we have for any $\varepsilon>0$ and $\theta=\max(\tau,4N^{-2})$ sufficiently small the estimate
\begin{equation}\label{h}
\sup\limits_{t^\prime\in[0,4\tau]}\Vert u^\theta _n\Vert_{X_\tau^{s_0,\frac12-\varepsilon}(T)}\leq C_T.
\end{equation}
\end{lemma}

\begin{proof}
By the estimate
$$
\Vert u_n\Vert_{X_\tau^{s,b}}\lesssim\Vert u\Vert_{X^{s,b}}+\tau^{b^\prime}\Vert u\Vert_{X^{s,1}},
$$
which can be found in \cite[Lemma 4.1]{Ji} or \cite[Lemma 3.4]{Rou}, with $b=\frac12-\varepsilon$ for $\varepsilon>0$ and $b^\prime=\frac12+\varepsilon$, we have
$$
\Vert u_n^\theta\Vert_{X_\tau^{s_0,b}} \lesssim\Vert u^\theta \Vert_{X^{s_0,b}}+\tau^{b^\prime}\Vert u^\theta \Vert_{X^{s_0,1}}.
$$
Note that we have the estimate (see \cite{Ji,JiO})
$$
\|u^\theta\|_{X^{s_0,1}(T)}\lesssim_{T,b^\prime}\theta^{-b^\prime}.
$$
Thus we get from \eqref{bmieux} that
$$
\Vert u_n^\theta\Vert_{X_\tau^{s_0,\frac12-\varepsilon}} \lesssim C_T,
$$
which ends the proof.
\end{proof}

In the following lemma, we shall give some important properties of $T_N$.

\begin{lemma}\label{tnprop}
For the operator $T_N$ introduced in Definition \ref{tndefine}, $\theta=\max(\tau,4N^{-2})$, and constants $b\in\mathbb{R}$, $c>0$, $p\in[1,\infty]$, we have the following estimates
\begin{align}
\label{tnkeep}T_N\Pi_\theta v_n&=\Pi_\theta v_n,\\
\label{lpparsv}\Vert T_N v_n\Vert_{l_\tau^pL^2}&=\Vert v_n(x_{jl})\Vert_{l_\tau^pl_h^2},\\
\label{discount}\Vert\Pi_{c\theta}v_n(x_{jl})\Vert_{l_\tau^pl_h^2}&\lesssim\Vert\Pi_{c\theta} v_n\Vert_{l_\tau^pL^2},\\
\label{discount2}\Vert T_N\Pi_{c\theta} v_n\Vert_{X_\tau^{0,b}}&\lesssim \Vert\Pi_{c\theta} v_n\Vert_{X_\tau^{0,b}}
\end{align}
for any appropriate function $v_n$.
\end{lemma}

\begin{proof}
By \eqref{tndef}, $T_N$ is a linear isometry (see also \cite{LiWu}), i.e. $T_N(u)$ is the unique function with frequencies $k_1,k_2$ satisfying $-\frac N2\leq k_1,k_2\leq\frac N2-1$ that coincides with $u$ at all spatial grid points $x_{jl}=(\tfrac{2\pi j}{N},\tfrac{2\pi l}{N}),~-\tfrac N2\leq j,l\leq\tfrac N2-1$.

We first observe that $\Pi_\theta v_n$ has no frequencies larger than $\theta^{-\frac12}\leq\frac N2$. Therefore, by the definition of $T_N$ (see Definition \ref{tndefine}), neither $\Pi_\theta v_n$ nor $T_N\Pi_\theta v_n$ has frequencies larger than $\frac N2$. Thus we can get \eqref{tnkeep} by the uniqueness.

As a direct corollary of \eqref{tnkeep}, we also point out that
\begin{equation}\label{iminuspi}
\Pi_\theta(T_Nv_n-v_n)=\Pi_\theta T_N(I-\Pi_\theta)v_n.
\end{equation}

To prove \eqref{lpparsv}, we first prove
\begin{equation}\label{fftmass}
\Vert T_N v_n\Vert_{L^2}=\Vert v_n(x_{jl})\Vert_{l_h^2},
\end{equation}
for every $n\in\mathbb{N}$, where $h=\frac1N$. Then \eqref{lpparsv} is a direct corollary by summing over time. By \eqref{parseval}, we have the spatial Parseval's identity
$$
\Vert T_N v_n\Vert_{L^2}=\Vert\mathcal{F}(T_N v_n)\Vert_{l^2},
$$
where $\mathcal{F}$ denotes the spatial Fourier transform, see \eqref{spaceft}. Note that, by \eqref{dft}, we have
$$
\mathcal{F}(T_N v_n)=h^2\mathcal{F}_N(T_N v_n)=h^2\mathcal{F}_N(v_n).
$$
Thus by \eqref{lpnorm}, it suffices to prove the discrete Parseval's identity
\begin{equation}\label{fftparsv}
\Vert\mathcal{F}_N(v_n)\Vert_{l^2}=N^2\Vert v_n(x_{jl})\Vert_{l^2_h}.
\end{equation}
If we define matrices $V=\Big(v_n(x_{jl})\Big)_{j,l=-\frac N2}^{\frac N2-1}$ and $G=\Big(\mathcal{F}_N(v_n)(k_1,k_2)\Big)_{k_1,k_2=-\frac N2}^{\frac N2-1}$, then by \eqref{dft}, we have
$$
G=PVP,
$$
where $P=(e^{-\frac{i2\pi jk}N})_{j,k=-\frac N2}^{\frac N2-1}$. Note that $P^*P=NI_N$, where $P^*$ denotes the adjoint matrix of $P$. Therefore we have
\begin{align*}
\Vert\mathcal{F}_N(v_n)\Vert_{l^2}^2&=\Vert G\Vert_F^2=\tr(G^*G)=\tr(P^*V^*P^*PVP)\\
&=N^2\tr(V^*V)=N^2\Vert V\Vert_F^2=N^2\Vert v_n(x_{jl})\Vert_{l^2}^2=N^4\Vert v_n(x_{jl})\Vert_{l^2_h}^2,
\end{align*}
where $\Vert\cdot\Vert_F$ denotes the Frobenius norm and $\tr$ denotes the trace. This concludes \eqref{fftparsv} and thus proves \eqref{lpparsv}.

To prove \eqref{discount}, we set $m$ as the smallest positive integer satisfying $(c\theta)^{-\frac12}\leq\frac{mN}2$, thus $c^{-\frac12}+1>m$. By \eqref{lpnorm}, \eqref{tnkeep} and \eqref{lpparsv}, we have
$$
\Vert\Pi_{c\theta}v_n(x_{jl})\Vert_{l_\tau^pl_h^2}\leq m\Vert\Pi_{c\theta}v_n(x_{ik})\Vert_{l_\tau^pl_{h/m}^2}=m\Vert T_{mN}\Pi_{c\theta}v_n(x_{ik})\Vert_{l_\tau^pl_{h/m}^2}=m\Vert\Pi_{c\theta}v_n\Vert_{l_\tau^pL^2},
$$
where $-\frac{mN}2\leq i,k\leq \frac{mN}2-1$.

To prove \eqref{discount2}, by \eqref{fftmass}, we have
$$
\Vert e^{-in\tau\Delta}T_N\Pi_{c\theta}v_n\Vert_{L^2}=\Vert e^{-in\tau\Delta}\Pi_{c\theta}v_n(x_{jl})\Vert_{l_h^2}.
$$
Summing over time, by \eqref{5}, we have for any $b\in\mathbb{R}$ that
$$
\Vert T_N\Pi_{c\theta}v_n\Vert_{X_\tau^{0,b}}\sim\Vert\langle D_\tau\rangle^b e^{-in\tau\Delta}T_N\Pi_{c\theta}v_n\Vert_{l_\tau^2L^2}=\Vert\langle D_\tau\rangle^b e^{-in\tau\Delta}\Pi_{c\theta}v_n(x_{jl})\Vert_{l_\tau^2l_h^2}.
$$
Again by \eqref{lpnorm}, \eqref{tnkeep} and \eqref{lpparsv}, we have
\begin{align*}
\Vert\langle D_\tau\rangle^b&e^{-in\tau\Delta}\Pi_{c\theta}v_n(x_{jl})\Vert_{l_\tau^2l_h^2}\leq m\Vert\langle D_\tau\rangle^b e^{-in\tau\Delta}\Pi_{c\theta}v_n(x_{ik})\Vert_{l_\tau^2l_{h/m}^2}\\
&=m\Vert\langle D_\tau\rangle^b e^{-in\tau\Delta}T_{mN}\Pi_{c\theta}v_n(x_{ik})\Vert_{l_\tau^2l_{h/m}^2}=m\Vert\langle D_\tau\rangle^b e^{-in\tau\Delta}\Pi_{c\theta}v_n\Vert_{l_\tau^2L^2}\sim\Vert\Pi_{c\theta}v_n\Vert_{X^{0,b}_\tau},
\end{align*}
where $-\frac{mN}2\leq i,k\leq \frac{mN}2-1$. This concludes \eqref{discount2}.
\end{proof}

Next, we shall give an estimate of $\mathcal{E}_{sloc}(t_k,N,\tau,u^\theta )$, see \eqref{localee}.

\begin{theorem}\label{theolocal}
For $s_0\leq2$ and $u^\theta$ as in Theorem~\ref{theou-utau}, we have for $\theta=\max(\tau,4N^{-2})$ sufficiently small that
$$
\Vert\mathcal{E}_{sloc}(t_k,N,\tau,u^\theta)\Vert_{X^{0,-b_1}_\tau}\leq C_T\tau\theta^{\frac{s_0}2},
$$
where $s_0$, $b_0$ fulfil the conditions in \eqref{s0b0}.
\end{theorem}

\begin{proof}
We first compute $\mathcal{E}_{sloc}$. By \eqref{tnkeep} and \eqref{iminuspi}, we have
\begin{equation}\label{sloc}
\begin{aligned}
\mathcal{E}_{sloc}(t_k,N,\tau,u^\theta)&=\Pi_\theta(T_N(e^{-i\tau|\Pi_\theta u^\theta(t_n)|^2}\Pi_\theta u^\theta(t_n))-(e^{-i\tau|\Pi_\theta  u^\theta(t_n)|^2}\Pi_\theta u^\theta(t_n)))\\
&=\Pi_\theta T_N\mathcal{E}_1(u^\theta(t_n))-\Pi_\theta\mathcal{E}_1(u^\theta(t_n))-i\tau\Pi_\theta T_N(I-\Pi_\theta)\mathcal{E}_2(u^\theta(t_n)),
\end{aligned}
\end{equation}
where
\begin{align*}
\mathcal{E}_1(u^\theta(t_n))&=(e^{-i\tau|\Pi_\theta u^\theta(t_n)|^2}-1+i\tau|\Pi_\theta u^\theta(t_n)|^2)\Pi_\theta u^\theta(t_n),\\
\mathcal{E}_2(u^\theta(t_n))&=|\Pi_\theta u^\theta(t_n)|^2\Pi_\theta u^\theta(t_n).
\end{align*}
Note that 
\begin{equation}\label{e2c9}
\Pi_{\frac\theta9}\mathcal{E}_2(u^\theta)=\mathcal{E}_2(u^\theta).
\end{equation}
By \eqref{b}, \eqref{ydual} and \eqref{lpparsv}, we get
\begin{equation}\label{te1prepare}
\Vert\Pi_\theta T_N\mathcal{E}_1(u^\theta(t_n))\Vert_{X_\tau^{0,-b_1}}\lesssim\tau^{-\frac12+b_1-\varepsilon}\Vert T_N\mathcal{E}_1(u^\theta(t_n))\Vert_{l_\tau^1L^2}=\tau^{-\frac12+b_1-\varepsilon}\Vert\mathcal{E}_1(u^\theta(t_n))\Vert_{l_\tau^1l_h^2},
\end{equation}
where $h=\frac1N$, and $\varepsilon>0$ is sufficiently small. Since $e^{\tau\alpha}-1-\tau\alpha\lesssim\tau^2|\alpha|^2$ for $\alpha\lesssim1$, we therefore have
\begin{equation}\label{te1start}
\Vert\mathcal{E}_1(u^\theta(t_n))\Vert_{l_\tau^1l_h^2}\lesssim\tau^2\Vert|\Pi_\theta u^\theta(t_n)|^4\Pi_\theta u^\theta(t_n)\Vert_{l_\tau^1l_h^2}=\tau^2\Vert(\Pi_\theta u^\theta(t_n))^5\Vert_{l_\tau^1l_h^2}.
\end{equation}
By \eqref{discount} and H\"older's inequality, we have
\begin{equation}\label{c12start}
\Vert(\Pi_\theta u^\theta(t_n))^5\Vert_{l_\tau^1l_h^2}\lesssim\Vert(\Pi_\theta u^\theta(t_n))^5\Vert_{l_\tau^1L^2}\lesssim\Vert \Pi_\theta u^\theta(t_n)\Vert^4_{l_\tau^4L^{20}}\Vert\Pi_\theta u^\theta(t_n)\Vert_{l_\tau^\infty L^\frac{10}3},
\end{equation}
since $\Pi_{\frac\theta{25}}(\Pi_\theta u^\theta(t_n))^5=(\Pi_\theta u^\theta(t_n))^5$.

For $s_0\in(0,\frac25]$, we get by Sobolev embeddings, \eqref{z} and \eqref{t} that
\begin{equation}\label{te1l4}
\Vert\Pi_\theta u^\theta(t_n)\Vert_{l_\tau^4L^{20}}\lesssim\Vert\Pi_\theta u^\theta(t_n)\Vert_{l_\tau^4W^{\frac25,4}}\lesssim\Vert\langle\partial_x \rangle^\frac25\Pi_\theta u^\theta(t_n)\Vert_{X_\tau^{5\varepsilon,\frac12-\varepsilon}}\lesssim\theta^{-\frac15-\frac{5\varepsilon}2+\frac{s_0}2}\Vert u^\theta(t_n)\Vert_{X_\tau^{s_0,\frac12-\varepsilon}}.
\end{equation}
Moreover, by Sobolev embedding, the continuity of $u^\theta$ and \eqref{d}, we get that
\begin{equation}\label{te1li}
\Vert\Pi_\theta u^\theta(t_n)\Vert_{l_\tau^\infty L^\frac{10}3}\lesssim\Vert\Pi_\theta u^\theta(t_n)\Vert_{l_\tau^\infty H^\frac25}\lesssim\theta^{-\frac15+\frac{s_0}2}\Vert u^\theta(t_n)\Vert_{L^\infty H^{s_0}}\lesssim\theta^{-\frac15+\frac{s_0}2}\Vert u^\theta(t_n)\Vert_{X^{s_0,\frac12+\varepsilon}}.
\end{equation}
Now, we combine \eqref{te1prepare}, \eqref{te1start}, \eqref{c12start}, \eqref{te1l4} and \eqref{te1li}, take $\varepsilon$ small enough ($s_0>50\varepsilon$, for example) and use \eqref{h} and \eqref{diffutau}, noting that $\tau\leq\theta$, to get
\begin{equation}\label{te1c1}
\Vert\Pi_\theta T_N\mathcal{E}_1(u^\theta(t_n))\Vert_{X^{0,-b_1}_\tau}\lesssim C_T\tau^{\frac32+b_1-\varepsilon}\theta^{\frac{5s_0}2-1-10\varepsilon}\lesssim C_T\tau\theta^{\frac{s_0}2}
\end{equation}
since $b_1>\frac12-\frac14s_0$.

For $s_0\in(\frac25,\frac32]$, we get by a similar argument as above that
\begin{equation}\label{te1c2}
\Vert\Pi_\theta T_N\mathcal{E}_1(u^\theta(t_n))\Vert_{X_\tau^{0,-b_1}}\lesssim\tau^{\frac32+b_1-\varepsilon}\Vert u^\theta(t_n)\Vert^4_{X_\tau^{s_0,\frac12-\varepsilon}}\Vert u^\theta(t_n)\Vert_{X^{s_0,\frac12+\varepsilon}}\lesssim C_T\tau\theta^{\frac34}\lesssim C_T\tau\theta^{\frac{s_0}2}
\end{equation}
since $b_1>\frac14$.

For $s_0\in(\frac32,2]$, by using similar techniques as \eqref{te1prepare} and \eqref{te1start}, H\"older's inequality and \eqref{lpparsv}, we obtain the rough estimate
$$
\Vert\Pi_\theta T_N\mathcal{E}_1(u^\theta(t_n))\Vert_{X_\tau^{0,-b_1}}\lesssim\tau^2\Vert(u^\theta(t_n))^5\Vert_{l_\tau^2L^2}\lesssim\tau^2\Vert u^\theta(t_n)\Vert^2_{l_\tau^4L^\infty}\Vert u^\theta(t_n)\Vert^3_{l_\tau^\infty L^6}.
$$
By Sobolev embedding and a similar argument as before, we have
\begin{equation}\label{te1c3}
\Vert\Pi_\theta T_N\mathcal{E}_1(u^\theta(t_n))\Vert_{X_\tau^{0,-b_1}}\lesssim\tau^2\Vert u^\theta(t_n)\Vert^2_{X_\tau^{s_0,\frac12-\varepsilon}}\Vert u^\theta(t_n)\Vert^3_{X^{s_0,\frac12+\varepsilon}}\lesssim C_T\tau^2\lesssim C_T\tau\theta^{\frac{s_0}2}.
\end{equation}
Combining \eqref{te1c1}, \eqref{te1c2} and \eqref{te1c3}, we get
\begin{equation}\label{te1}
\Vert\Pi_\theta T_N\mathcal{E}_1(u^\theta(t_n))\Vert_{X^{0,-b_1}_\tau}\lesssim C_T\tau\theta^{\frac{s_0}2}.
\end{equation}

Similarly, we get that
\begin{equation}\label{e1}
\Vert\Pi_\theta\mathcal{E}_1(u^\theta(t_n))\Vert_{X_\tau^{0,-b_1}}\leq C_T\tau\theta^{\frac{s_0}2}.
\end{equation}

Finally, by \eqref{discount2}, \eqref{e2c9} and \eqref{t}, we have
$$
\Vert\Pi_\theta T_N(I-\Pi_\theta)\mathcal{E}_2(u^\theta(t_n))\Vert_{X_\tau^{0,-b_1}}\lesssim \Vert (I-\Pi_\theta)\mathcal{E}_2(u^\theta(t_n))\Vert_{X_\tau^{0,-b_1}}\lesssim \theta^{\frac{s_0}2}\Vert\mathcal{E}_2(u^\theta(t_n))\Vert_{X_\tau^{s_0,-b_1}}.
$$
Thus by \eqref{r9} and \eqref{h}, we get that
\begin{equation}\label{e2}
\begin{aligned}
\Vert i\tau\Pi_\theta T_N(I-\Pi_\theta)\mathcal{E}_2(u^\theta(t_n))\Vert_{X_\tau^{0,-b_1}}&\lesssim\tau\theta^{\frac{s_0}2}\Vert|\Pi_\theta u^\theta(t_n)|^2\Pi_\theta u^\theta(t_n)\Vert_{X_\tau^{s_0,-b_1}}
\\&\lesssim\tau\theta^{\frac{s_0}2}\Vert u^\theta(t_n)\Vert_{X_\tau^{s_0,b_1}}^3\leq C_T\tau\theta^{\frac{s_0}2}.
\end{aligned}
\end{equation}
The theorem concludes by combining \eqref{sloc}, \eqref{te1}, \eqref{e1}, and \eqref{e2}.
\end{proof}

\section{Proof of Theorem \ref{mainthm}}\label{sectionglobal}

In this section, we will first estimate the global error $e_n$. Similar to \cite[sect.~3]{Ost1}, we write the global error as
\begin{equation}\label{k}
\begin{aligned}
e_n&=u^\theta(t_n)-u_n\\
&=u^\theta(t_{n-1})-u_{n-1}-i\tau e^{i\tau\Delta}(\Phi^\tau_N(u^\theta(t_{n-1}))-\Phi^\tau_N(u_{n-1}))-ie^{i\tau\Delta}\mathcal{E}_{loc}(t_{n-1},N,\tau,u^\theta)\\
&=-i\tau\sum\limits_{k=0}^{n-1}e^{i(n-k)\tau\Delta}(\Phi^\tau_N(u^\theta(t_k))-\Phi^\tau_N(u_k))-i\sum\limits_{k=0}^{n-1}e^{i(n-k)\tau\Delta}\mathcal{E}_{loc}(t_k,N,\tau,u^\theta )
\end{aligned}
\end{equation}

with the nonlinear flow
\begin{equation}\label{defphitaun}
\Phi_N^\tau(w)=-\Pi_\theta T_N\Big(\dfrac{e^{-i\tau|\Pi_\theta w|^2}-1}{i\tau}\Pi_\theta w\Big).
\end{equation}

\begin{proposition}\label{propen}
For $s_0\leq2$ and $u^\theta$ as in Theorem~\ref{theou-utau}, we have for $\theta=\max(\tau,4N^{-2})$ sufficiently small that
$$
\Vert e_n\Vert_{X^{0,b_0}_{\tau}(T)} \lesssim\theta^\frac{s_0}2,
$$
where $s_0$, $b_0$ fulfil the conditions in \eqref{s0b0}.
\end{proposition}

\begin{proof}
Take a smooth function $\eta$ which is 1 on $[-1,1]$ and compactly supported in $[-2,2]$. In the proof we shall still denote the error of the truncated version of \eqref{k} by $e_n$, given as
\begin{equation}\label{m}
e_n=-i\tau\eta(t_n)\sum\limits_{k=0}^{n-1}e^{i(n-k)\tau\Delta}\,\eta\Big(\frac{t_k}{T_1}\Big)\big(\Phi_N^\tau(u^\theta(t_k))-\Phi_N^\tau(u^\theta(t_k)-e_k)\big)+\mathcal{R}_n,
\end{equation}
where
$$
\mathcal{R}_n=-i\eta(t_n)\sum\limits_{k=0}^{n-1}e^{i(n-k)\tau\Delta}\eta\Big(\frac{t_k}{T_1}\Big)\mathcal{E}_{loc}(t_k,N,\tau,u^\theta).
$$
Note that for $0\leq n\leq N_1$, where $N_1=[\frac{T_1}{\tau}]$ with $T_1\leq\min(1,T)$ to be chosen, this indeed coincides with $u^\theta(t_n)- u_n$. By using \eqref{q}, \eqref{timeloc} and Theorem~\ref{theolocal}, we get
\begin{equation}\label{8}
\Vert\mathcal{R}_n\Vert_{X^{0,b_0}_\tau}\lesssim\tau^{-1}\Vert\mathcal{E}_{loc}(t_n,N,\tau,u^\theta)\Vert_{X^{0,b_0-1}_\tau}\lesssim C_T\theta^{\frac{s_0}2}.
\end{equation}

Thus by \eqref{m} and \eqref{8}, we have
$$
\Vert e_n\Vert_{X_\tau^{0,b_0}}\lesssim\left\Vert\tau\eta(t_n)\sum\limits_{k=0}^{n-1}e^{i(n-k)\tau\Delta}\eta\Big(\frac{t_k}{T_1}\Big)\big(\Phi^\tau_N(u^\theta(t_k))-\Phi^\tau_N(u^\theta(t_k)-e_k)\big)\right\Vert_{X_\tau^{0,b_0}}+C_T\theta^{\frac{s_0}2}.
$$
By the definition of the nonlinear flow, see \eqref{defphitaun}, we can write $\Phi_N^\tau$ in the form
$$
\Phi_N^\tau(w)=\Pi_\theta T_N F_N^\tau(w) + \Pi_\theta T_N R_N^\tau(w),
$$
where
\begin{align*}
F_N^\tau(w)&= |\Pi_\theta(w)|^2\Pi_\theta w, \\
R_N^\tau(w)& = F_N^\tau(w) R_N^{\tau,1}(w), \quad R_N^{\tau,1}(w)= \int_0^1 \left(e^{-is\tau |\Pi_\theta w|^2}-1\right) ds.
\end{align*}
Note that $R_N^\tau(w)$ behaves like a quintic term in the sense that
$$
|R_N^\tau(w)| \lesssim \tau|\Pi_\theta w|^5
$$
with the gain of a factor $\tau$.

By using this decomposition and \eqref{9}, \eqref{q}, we get that
\begin{equation}\label{esten1}
\begin{aligned}
\Vert e_n\Vert_{X^{0,b_0}_\tau}\leq C_TT_1^{\varepsilon_0}\Vert \Pi_\theta&T_N\bigl(F_N^\tau(u^\theta(t_n))-F_N^\tau(u^\theta(t_n)-e_n)\bigr)\Vert_{X_\tau^{0,-b^\prime}} \\
&+ C_T\Vert\Pi_\theta T_N\bigl(R_N^\tau(u^\theta(t_n))-R^\tau_N(u^\theta(t_n)-e_n) \bigr)\Vert_{X_\tau^{0,-b_1}}+C_T\tau^{\frac{s_0}2},
\end{aligned}
\end{equation}
where $b^\prime$ is chosen such that $\max(\frac14,\frac12-\frac14s_0)<b^\prime<b_1$, $\varepsilon_0=b_1-b^\prime$.

Now we expand
$$
F_N^\tau(u^\theta(t_n))-F_N^\tau(u^\theta(t_n)-e_n)= L^\theta_n + Q^\theta_n + C^\theta_n
$$
with
\begin{align*}
L^\theta_n&=2|\Pi_\theta u^\theta(t_n)|^2\Pi_\theta e_n +\Pi_\theta u^\theta(t_n)\Pi_\theta u^\theta(t_n)\Pi_\theta\overline{e}_n, \\
Q^\theta_n&=-\Pi_\theta\overline{u}^\theta(t_n)\Pi_\theta e_n \Pi_\theta e_n-2\Pi_\theta u^\theta(t_n)|\Pi_\theta e_n|^2, \\
C^\theta_n&=|\Pi_\theta e_n|^2 \Pi_\theta e_n.
\end{align*}
Note that we have $L^\theta_n=\Pi_{\frac\theta9}L^\theta_n,~Q^\theta_n=\Pi_{\frac\theta9}Q^\theta_n$ and $C^\theta_n=\Pi_{\frac\theta9}C^\theta_n$. Thus we get by using \eqref{discount2} and \eqref{rbis9} that
\begin{equation}\label{ouf4}
\begin{aligned}
\Vert\Pi_\theta T_N\bigl(F_N^\tau(u^\theta(t_n))&-F_N^\tau(u^\theta(t_n)-e_n)\bigr)\Vert_{X_\tau^{0,-b^\prime}}\lesssim\|u^\theta(t_n)\|_{X^{s_0,b_1}_\tau}^2\|e_n\|_{X^{0,b_1}_\tau}\\
&+\|u^\theta(t_n)\|_{X^{s_0,b_1}_\tau}\|e_n\|_{X^{0,b_1}_\tau} \|\Pi_\theta e_n\|_{X^{s_1,b_1}_\tau}+\|\Pi_\theta e_n\|^2_{X^{s_1,b_1}_\tau}\|e_n\|_{X^{0,b_1}_\tau},
\end{aligned}
\end{equation}
with $s_1\in (0,s_0]$ to be chosen. Note that we have freely added $\Pi_\theta$ in the second and the third term of the right-hand side, which is allowed since $\Pi_\theta e_n=e_n$. By using \eqref{t}, which reads
$$
\|\Pi_\theta e_n\|_{X^{s_1,b_1}_\tau}\lesssim\theta^{-\frac{s_1}2}\|e_n\|_{X^{0,b_1}_\tau},
$$
and \eqref{h}, we thus get that
\begin{equation}\label{esten2}
\begin{aligned}
\Vert\Pi_\theta T_N\bigl(F_N^\tau(u^\theta(t_n))-F_N^\tau(&u^\theta(t_n)-e_n)\bigr)\Vert_{X_\tau^{0,-b^\prime}} \\
& \leq C_T \|e_n\|_{X^{0,b_1}_\tau}+C_{T,s_1}\theta^{-\frac{s_1}2}\|e_n\|_{X^{0,b_1}_\tau}^2+C_{T,s_1}\theta^{-s_1}\|e_n\|_{X^{0,b_1}_\tau}^3.
\end{aligned}
\end{equation}
We shall now estimate $\Vert\Pi_\theta T_N\left(R_N^\tau(u^\theta(t_n))-R_N^\tau(u^\theta(t_n)-e_n)\right)\Vert_{X_\tau^{0,-b_1}}$. We first use \eqref{b} and \eqref{ydual} to get
\begin{equation}\label{rnouf}
\begin{aligned}
\Vert\Pi_\theta T_N\bigl(R_N^\tau(u^\theta(t_n))&-R_N^\tau(u^\theta(t_n)-e_n)\bigr)\Vert_{X_\tau^{0,-b_1}}\\
&\lesssim\tau^{-\delta-(\frac12-b_1)}\Vert\Pi_\theta T_N\bigl(R_N^\tau (u^\theta(t_n))-R_N^\tau(u^\theta(t_n)-e_n)\bigr)\Vert_{X_\tau^{0,-\frac12-\delta}}\\
&\lesssim\tau^{-\delta-(\frac12-b_1)}\Vert T_N\big(R_N^\tau(u^\theta(t_n))-R_N^\tau(u^\theta(t_n)-e_n)\big)\Vert_{l^1_\tau L^2},
\end{aligned}
\end{equation}
where $\delta>0$ will be chosen sufficiently small. Then, since we have the pointwise (uniform in $\tau$) estimate
$$
\left|R_N^\tau(u^\theta(t_n))-R_N^\tau(u^\theta(t_n)-e_n)\right|\lesssim \tau\sum_{j=1}^5 |\Pi_\theta u^\theta(t_n)|^{5-j} |\Pi_\theta e_n|^j,
$$
we thus deduce by \eqref{rnouf}, \eqref{lpparsv} and \eqref{discount} that
\begin{equation}\label{ouf3}
\begin{aligned}
\Big\Vert\Pi_\theta T_N\big(R_N^\tau(u^\theta(t_n))-R_N^\tau(u^\theta(t_n)-e_n)\big)\Big\Vert_{X_\tau^{0,-b_1}}&\lesssim\tau^{-\delta-(\frac12-b_1)}\Vert R_N^\tau(u^\theta(t_n))-R_N^\tau(u^\theta(t_n)-e_n)\Vert_{l^1_\tau l_h^2}\\
&\lesssim\tau^{1-\delta-(\frac12-b_1)}\sum_{j=1}^5 \big\||\Pi_\theta u^\theta(t_n)|^{5-j}|\Pi_\theta e_n|^j\big\|_{l^1_\tau l_h^2}\\
&\lesssim \theta^{1-\delta-(\frac12-b_1)}\sum_{j=1}^5 \big\| |\Pi_\theta u^\theta(t_n)|^{5-j} |\Pi_\theta e_n|^j \big\|_{l^1_\tau L^2}.
\end{aligned}
\end{equation}
By H\"{o}lder's inequality, we get that
\begin{align*}
\sum_{j=1}^5 \big\| |\Pi_\theta u^\theta(&t_n)|^{5-j} |\Pi_\theta e_n|^j \big\|_{l^1_\tau L^2} \lesssim
\|\Pi_\theta u^\theta(t_n)\|_{l^4_\tau L^\infty}^4\|\Pi_\theta e_n\|_{l^\infty_\tau L^2}\\[-2.5mm]
&+ \|\Pi_\theta u^\theta(t_n)\|_{l^4_\tau L^\infty}^3\|\Pi_\theta e_n\|_{l^4_\tau L^\infty}\|\Pi_\theta e_n\|_{l^\infty_\tau L^2}
+ \|\Pi_\theta u^\theta(t_n)\|_{l^4_\tau L^\infty}^2\|\Pi_\theta e_n\|_{l^4_\tau L^\infty}^2\|\Pi_\theta e_n\|_{l^\infty_\tau L^2}\\
& +\|\Pi_\theta u^\theta(t_n)\|_{l^4_\tau L^\infty}\|\Pi_\theta e_n\|_{l^4_\tau L^\infty}^3\|\Pi_\theta e_n\|_{l^\infty_\tau L^2} +  \|\Pi_\theta e_n\|_{l^4_\tau L^\infty}^4\|\Pi_\theta e_n\|_{l^\infty_\tau L^2}.
\end{align*}
To estimate the right-hand side, we first get by \eqref{y} that
$$
\|\Pi_\theta e_n\|_{l^\infty_\tau L^2}\lesssim\|e_n\|_{X^{0,b_0}_\tau}.
$$
Next, by Sobolev embeddings and \eqref{z}, \eqref{t}, we get that
$$
\|\Pi_\theta e_n\|_{l^4_\tau L^\infty}\lesssim \|\Pi_\theta e_n\|_{l^4_\tau W^{\frac12+\delta,4}}\lesssim \|\Pi_\theta e_n\|_{X_\tau^{\frac12+2\delta,b_0}}\lesssim \theta^{-\frac14-\delta}\|e_n\|_{X_\tau^{0,b_0}}.
$$
Finally, by a similar argument, we get that
$$
\|\Pi_\theta u^\theta(t_n)\|_{l^4_\tau L^\infty} \lesssim \|\Pi_\theta u^\theta(t_n)\|_{l^4_\tau W^{\frac12+\delta,4}} \lesssim \|u^\theta(t_n)\|_{X^{\frac12+6\delta,\frac12-\delta}_\tau}\lesssim \theta^{\frac{s_0}2-3\delta-\frac14} C_T
$$
for $s_0\in(0,\frac12]$. In the case $s_0>\frac12$, we just have a simpler estimate directly from the Sobolev embedding. To be specific, by using a similar argument as before, we get
$$
\|\Pi_\theta u^\theta(t_n)\|_{l^4_\tau L^\infty} \lesssim \|\Pi_\theta u^\theta(t_n)\|_{l^4_\tau W^{\frac12+\delta,4}} \lesssim \|u^\theta(t_n)\|_{X^{s_0,\frac12-\delta}_\tau}\lesssim C_T
$$
by taking $\delta<\frac15(s_0-\frac12)$.
This yields in the case $s_0\leq\frac12$ that
\begin{equation}\label{ouf}
\begin{aligned}
\sum_{j=1}^5 \left\| |\Pi_\theta u^\theta(t_n)|^{5-j} |\Pi_\theta e_n|^j \right\|_{l^1_\tau L^2} &\lesssim C_{T,\delta}\left(\theta^{2s_0-1-12\delta}\|e_n\|_{X_\tau^{0,b_0}}+\theta^{\frac32s_0-1-10\delta} \|e_n\|_{X_\tau^{0,b_0}}^2 \right. \\
&\hspace{-5mm}\left. +\theta^{s_0-1-8\delta} \|e_n\|_{X_\tau^{0,b_0}}^3+ \theta^{\frac{s_0}2-1-6\delta}\|e_n\|_{X_\tau^{0,b_0}}^4 +\theta^{-1-4\delta}\|e_n\|_{X_{\tau}^{0,b_0}}^5\right).
\end{aligned}
\end{equation}
The case $s_0>\frac12$ can be handled by similar arguments. Consequently, by combining \eqref{esten1}, \eqref{esten2}, \eqref{ouf3}, and \eqref{ouf}, we obtain that
\begin{multline*}
\|e_n\|_{X^{0,b_0}_\tau} \leq C_T \left(( T_1^{\varepsilon_0} + C_{T, \delta}\theta^{2s_0-13\delta-(\frac12-b_1)})\|e_n\|_{X^{0,b_0}_\tau} \right. \\
\left.{} + C_{T, s_1,\delta} (\theta^{-\frac{s_1}2} + \theta^{\frac32s_0-11\delta-(\frac12-b_1)}) \|e_n\|_{X^{0, b_0}_\tau}^2 + C_{T, s_1,\delta}(\theta^{-s_1}+\theta^{s_0-9\delta-(\frac12-b_1)}) \|e_n\|_{X^{0, b_0}_\tau}^3 \right. \\
\left.{} + C_{T,\delta} \theta^{\frac{s_0}2-7\delta-(\frac12-b_1)} \|e_n\|_{X^{0,b_0}_\tau}^4+ C_{T,\delta} \theta^{-5\delta-(\frac12-b_1)}\|e_n\|_{X^{0, b_0}_\tau}^5+\theta^{\frac{s_0}2}\right).
\end{multline*}
To conclude, we observe that for $s_0>0$ fixed, we first choose $s_1$ satisfying $\max(\tfrac14,\tfrac12-\tfrac14s_1)<b^\prime<b_1$ and $\tfrac{s_0}2-\tfrac{s_1}2>0$, i.e. $2-4b_1<2-4b^\prime<s_1<s_0$ (this is possible since we always have $b^\prime>\tfrac12-\tfrac14s_0$), and then $\delta$ sufficiently small ($s_1>50\delta$ for example). Since $\tfrac12-b_1<\tfrac14s_1$, we get
\begin{multline*}
\|e_n\|_{X^{0,b_0}_\tau} \leq C_T \left(( T_1^{\varepsilon_0} + C_{T,s_1, \delta} \theta^\rho) \|e_n\|_{X^{0, b_0}_\tau} + C_{T,s_1,\delta}\theta^{-\frac{s_1}2} \|e_n\|_{X^{0, b_0}_\tau}^2 + C_{T,s_1,\delta}\theta^{-s_1} \|e_n\|_{X^{0, b_0}_\tau}^3 \right.\\
\left.+ C_{T,\delta} \theta^{-\frac32s_1}\|e_n\|_{X^{0,b_0}_\tau}^4+ C_{T,\delta}\theta^{-2s_1}\|e_n\|_{X^{0,b_0}_\tau}^5+\tau^\frac{s_0}2\right)
\end{multline*}
for some $\rho>0$. We then choose $T_1$, $\tau$ and $4N^{-2}$ small enough with respect to $C_T$, so that $C_T( T_1^{\varepsilon_0} + C_{T,s_1,\delta} \theta^\rho)\leq\frac12$.
This yields
\begin{align*}
\|e_n\|_{X^{0, b_0}_\tau} \leq C_{T,s_1,\delta}\theta^{-\frac{s_1}2} &\|e_n\|_{X^{0, b_0}_\tau}^2 +C_{T,s_1,\delta} \theta^{-s_1} \|e_n\|_{X^{0, b_0}_\tau}^3 \\
&+C_{T,\delta}\theta^{-\frac32s_1}\|e_n\|_{X^{0, b_0}_\tau}^4+C_{T,\delta}\theta^{-2s_1}\|e_n\|_{X^{0,b_0}_\tau}^5+C_T\theta^\frac{s_0}2.
\end{align*}
We next deduce that for $\tau$ and $4N^{-2}$ sufficiently small
$$
\|e_n\|_{X^{0,b_0}_\tau} \leq C_T \theta^\frac{s_0}2.
$$
This proves the desired estimate for $0 \leq n\tau \leq T_1$. Since the choice of $T_1$ depends only on $T$, we can then reiterate the argument on $T_1\leq n\tau\leq2T_1$ and so on, to get finally the estimate for $0\leq n\tau\leq T$ for $\tau$ sufficiently small.
\end{proof}

We finally prove our main result, Theorem \ref{mainthm}.

\begin{proof}
To estimate the error $\Vert u(t_n)-u_n\Vert_{L^2}$, we just use that
$$
\Vert u(t_n)-u_n\Vert_{L^2}\leq\Vert u(t_n)-u^\theta(t_n)\Vert_{L^2}+\Vert u^\theta(t_n)-u_n\Vert_{L^2}\leq\Vert u-u^\theta \Vert_{L^\infty([0, T],  L^2)}+\Vert e_n\Vert_{l^\infty_{\tau}(0 \leq n\tau \leq T, L^2)}
$$
Next, we use the embeddings \eqref{d} and \eqref{y} combined  with Theorem \ref{theou-utau}, Proposition \ref{propen} to get that
$$
\Vert u(t_n)-u_n\Vert_{L^2} \leq C_T\theta^{\frac{s_0}2}.
$$
This concludes the proof of Theorem~\ref{mainthm}.
\end{proof}

\section{Numerical experiments}\label{sectionnumerexp}

\begin{figure}[h]
\begin{center}
\subfigure[]{\includegraphics[width=0.49\textwidth]{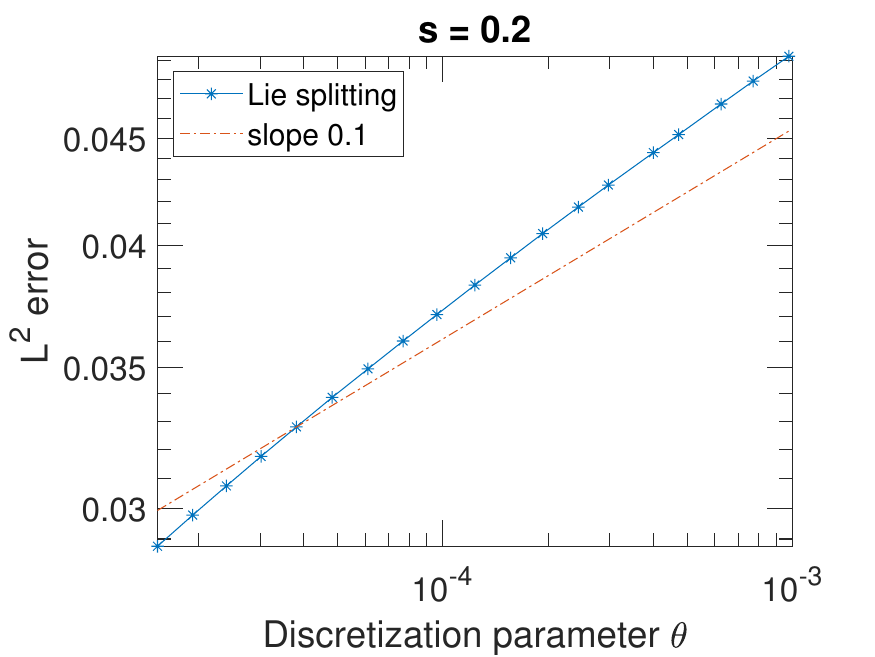}}
\subfigure[]{\includegraphics[width=0.49\textwidth]{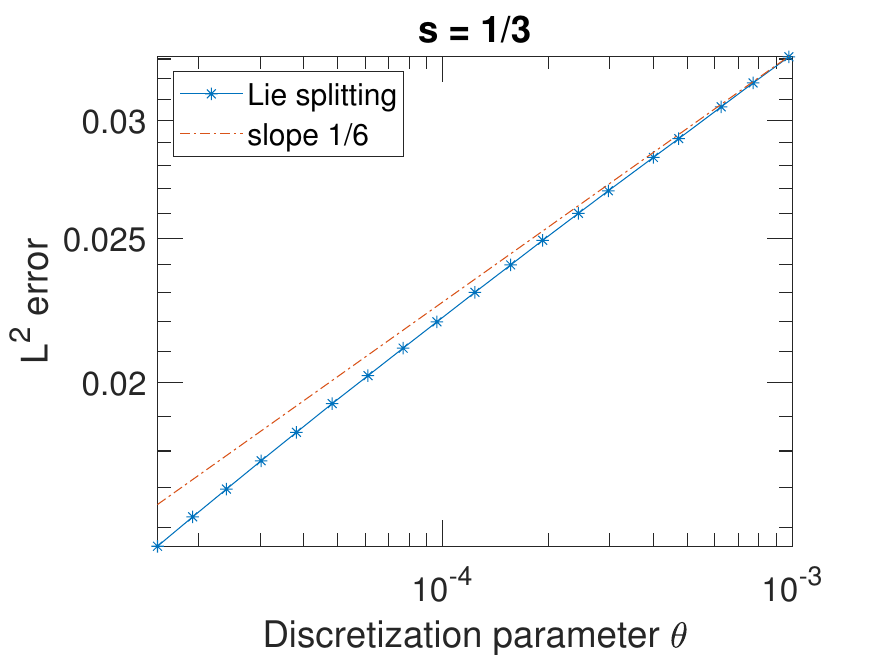}}
\smallskip
\subfigure[]{\includegraphics[width=0.49\textwidth]{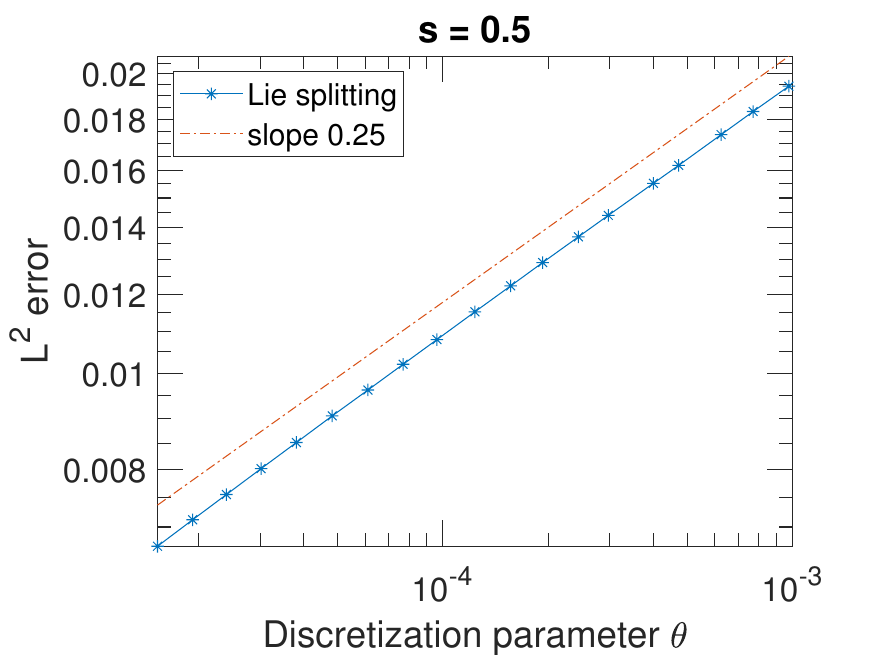}}
\subfigure[]{\includegraphics[width=0.49\textwidth]{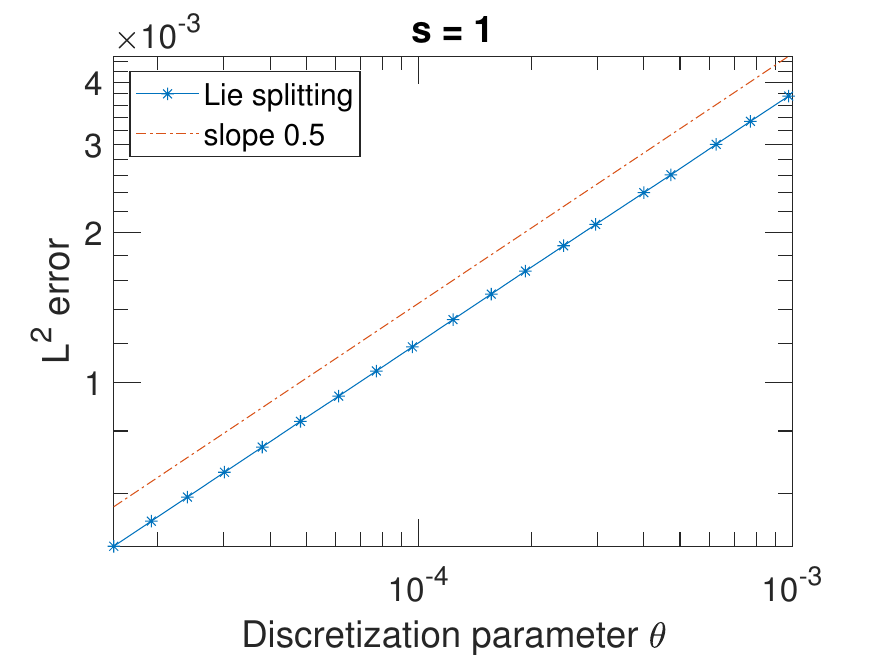}}
\end{center}
\caption{$L^2$ error of the fully discrete Lie splitting scheme \eqref{liespace} for rough initial data $u_0\in H^s(\mathbb{T}^2)$. We took the choice $\theta=\tau=4N^{-2}$. \newline
(a) $s=0.2$; \quad (b) $s=1/3$; \quad (c) $s=0.5$; \quad (d) $s=1$.\label{figure1}}
\end{figure}

In this section, we numerically illustrate our main result (Theorem~\ref{mainthm}). We display the convergence order of the fully discrete filtered Lie splitting method \eqref{liespace} with rough initial data. For this purpose, we consider the periodic NLS (\ref{nls}) in two dimensions, discretized with the filtered Lie splitting method \eqref{liespace} and initial data
$$
u_0=\sum\limits_{k\in\mathbb{Z}^2}\langle k\rangle^{-(s+1+\varepsilon)}\tilde{g}_ke^{i\langle k,x\rangle}\in H^s,
$$
where $\tilde{g}_k$ are random variables uniformly distributed in the square $[-1,1]+i[-1,1]$ and $\varepsilon>0$ arbitrarily small. In our experiments, we set $\varepsilon=0$. Moreover, we normalize $\|u_0\|_{L^2}=0.1$.

In Figure~\ref{figure1}, we employ a standard Fourier pseudospectral method for the discretization in space. For the reference solution, we choose as largest Fourier mode $K=(2^{13},2^{13})$, i.e., the spatial mesh size $\Delta x=0.0008$ for both of the spatial directions, and we use the Lie splitting method with $K$ spatial points and a very small time step size $\tau=2^{-24}$. We choose $T=0.25$ to be the final time. For the numerical tests, we choose $\tau$ between $2^{-16}$ and $2^{-10}$, and keep $N=2\tau^{-\frac12}$, i.e. $\theta=\tau=4N^{-2}$.

Figure~\ref{figure1} clearly shows that our numerical experiments confirm the convergence rate of order $\mathcal{O}(\theta^{\frac s2})$ for solutions in $H^s$ (see Theorem~\ref{mainthm}) with $s=0.2,~1/3,~0.5$ and 1.

\begin{figure}[t]
\begin{center}
\subfigure[]{\includegraphics[width=0.49\textwidth]{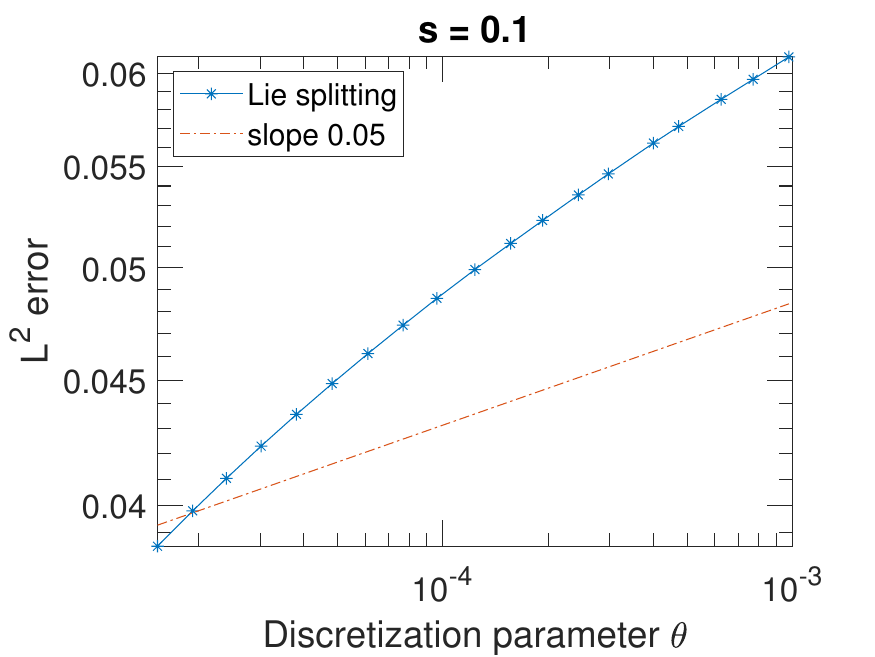}}
\subfigure[]{\includegraphics[width=0.49\textwidth]{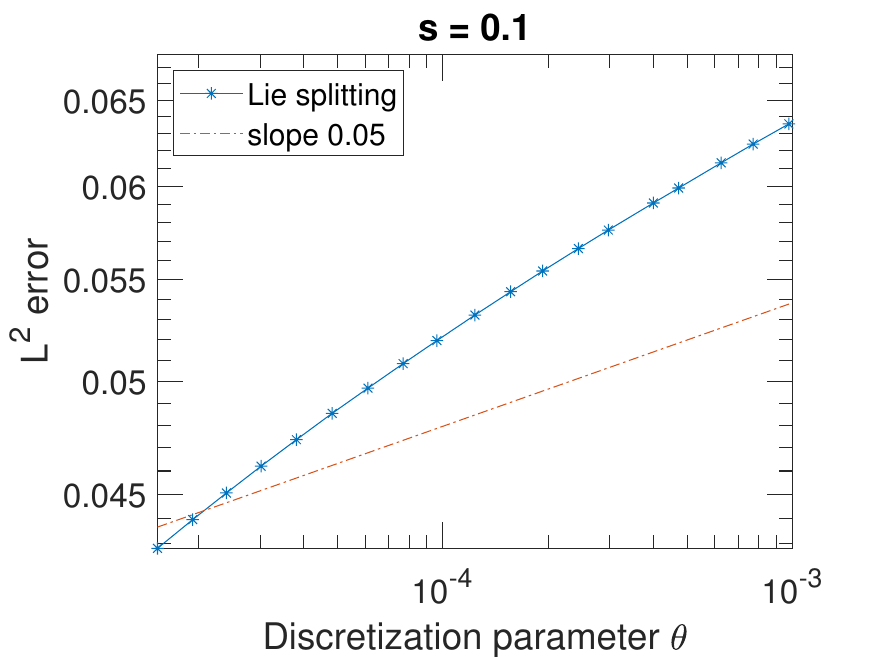}}
\end{center}
\caption{$L^2$ error of the fully discrete Lie splitting scheme for rough initial data $u_0\in H^{0.1}$ with different reference solutions. We took the choice $\theta=\tau=4N^{-2}$. (a) Reference solution with largest Fourier mode $K=(2^{12},2^{12})$, spatial mesh size $\Delta x=0.0015$, and time step size $\tau=2^{-22}$; (b) Reference solution with largest Fourier mode $K=(2^{13}, 2^{13})$, spatial mesh size $\Delta x=0.0008$, and time step size $\tau=2^{-24}$.\label{figure2}}
\end{figure}

We also did some experiments for very small $s$. As an example, we discuss here the case $s=0.1$. The convergence in this case is very slow. In order to obtain an accurate reference solution for small $s$ would require an enormous number of Fourier modes. In Figure~\ref{figure2}, where we have used the same scales in both subfigures, we demonstrate that the numerical result approaches the proven order of convergence as the number of Fourier modes increases. The proven order, however, only shows up for sufficiently high spatial resolution, which goes beyond the capacity of the used GPU cluster.

\section{Acknowledgement}

KS gratefully acknowledges funding from the European Research Council (ERC) under the European Union's Horizon 2020 research innovation programme (grant agreement No. 850941).

{}
\end{document}